\documentclass[12pt,reqno]{amsart}
\usepackage{amsmath, amsfonts, amssymb, amsthm, hyperref}
\usepackage{bm}
\allowdisplaybreaks[4]
\def\l{\left}
\def\r{\right}
\def\bg{\bigg}
\def\({\bg(}
\def\){\bg)}
\def\t{\text}
\def\f{\frac}

\def\ls{\leqslant}
\def\gs{\geqslant}

\def\bi{\binom}
\def\al{\alpha}

\def\eq{\equiv}

\def\Z{\mathbb Z}

\def\N{\mathbb N}

\def\1{{\bf 1}}

\def\pmod #1{\ ({\rm{mod}}\ #1)}

\def\<{\langle}
\def\>{\rangle}
\theoremstyle{plain}
\newtheorem{theorem}{Theorem}[section]
\newtheorem{lemma}{Lemma}

\newtheorem{conjecture}{Conjecture}

\theoremstyle{definition}

\theoremstyle{remark}
\newtheorem{remark}{Remark}

\makeatletter
\@namedef{subjclassname@2020}{%
  \textup{2020} Mathematics Subject Classification}
\makeatother
 \vspace{4mm}

\begin{document}
\hbox{Preprint}
\medskip

\title[On $p$-adic congruences involving $\sqrt d$]
{On $p$-adic congruences involving $\sqrt d$}
\author{Bo Jiang}
\address {(Bo Jiang) Department of Mathematics, Nanjing
University, Nanjing 210093, People's Republic of China}
\email{211501013@smail.nju.edu.cn}

\author{Zhi-Wei Sun}
\address {(Zhi-Wei Sun, corresponding author) Department of Mathematics, Nanjing
University, Nanjing 210093, People's Republic of China}
\email{zwsun@nju.edu.cn}

\keywords{$p$-adic congruences, polynomials over finite fields, quadratic fields.
\newline \indent Supported by the National Natural Science Foundation of China (grant no. 12371004).}
\begin{abstract}
Let $p$ be an odd prime and let $d$ be an integer not divisible by $p$.
We prove that
$$ \prod_{1\ls m,n\ls p-1\atop p\nmid m^2-dn^2}(x-(m+n\sqrt{d}))
\eq \begin{cases}\sum_{k=1}^{p-2}\frac{k(k+1)}2x^{(k-1)(p-1)}\pmod p
&\text{if}\ (\frac dp)=1,\\\sum_{k=0}^{(p-1)/2}x^{2k(p-1)}
\pmod p&\text {if}\ (\frac dp)=-1,
\end{cases}$$ 
where $(\frac dp)$ denotes the Legendre symbol. This extends a recent conjecture of N. Kalinin.
We also obtain the Wolstenholme-type congruence
$$\sum_{1\ls m,n\ls p-1\atop p\nmid m^2-dn^2}\frac1{m+n\sqrt d}\equiv0\pmod{p^2}.$$

\end{abstract}
\maketitle

\section{Introduction}
\setcounter{lemma}{0}
\setcounter{theorem}{0}
\setcounter{equation}{0}
\setcounter{conjecture}{0}
\setcounter{remark}{0}
\setcounter{corollary}{0}

In 1862, J. Wolstenholme \cite{W} proved for any prime $p>3$ the classical congruence
$$\sum_{n=1}^{p-1}\f1 n\eq0\pmod {p^2}.$$
Quite recently, N. Kalinin \cite{K} extended Wolstenholme's congruence to the ring
$$\Z[i]=\{m+ni:\ m,n\in\Z\}$$
of Gaussian integers. Namely, he deduced that for any prime $p>5$ we have
$$\sum_{1\ls m,n\ls p-1\atop p\nmid m^2+n^2}\f1{m+ni}\eq0\pmod{p^4}.$$
In the same paper, N. Kalinin made the following conjecture:

\begin{conjecture} Let $p>3$ be a prime, and set
$$g_p(x):=\prod_{1\ls m,n\ls p-1\atop p\nmid m^2+n^2}(x-(m+ni)).$$
If $p\eq1\pmod4$, then
$$g_p(x)\eq 1+\sum_{k=1}^{p-3}c_kx^{k(p-1)}\pmod p$$
for some $c_1,\ldots c_{p-3}\in\Z$ with $c_{p-3}=1$.
If $p\eq 3\pmod 4$, then
$$g_p(x)\eq  1+x^{2(p-1)}+x^{4(p-1)}+\dots+x^{(p-1)^2}=\frac{1-x^{p^2-1}}{1-x^{2(p-1)}}\pmod p.$$
\end{conjecture}

In this paper we establish the following result which extends the above conjecture of Kalinin.

\begin{theorem}\label{Th1.1}
	Let $p$ be an odd prime and let $d\in\Z$ with $p\nmid d$. For the polynomial
\begin{equation}\label{P}P(x):=\prod_{1\ls m,n\ls p-1\atop p\nmid m^2-dn^2}(x-(m+n\sqrt{d})),
\end{equation}
we have
\begin{equation}\label{P-cong}P(x)\eq \begin{cases}(x^{p-1}-1)^{p-3}\eq\sum_{k=1}^{p-2}\f{k(k+1)}2x^{(k-1)(p-1)}\pmod p
&\t{if}\ (\f dp)=1,\\(x^{p^2-1}-1)/(x^{2(p-1)}-1)=\sum_{k=0}^{(p-1)/2}x^{2k(p-1)}
\pmod p&\t{if}\ \f dp)=-1.
\end{cases}\end{equation}
\end{theorem}

We also obtain the following Wolstenholme-type result.

\begin{theorem}\label{Th1.2} For any odd prime $p$ and any integer $d\not\eq0\pmod p$, we have
\begin{equation}\label{p2}\sum_{1\ls m,n\ls p-1\atop p\nmid m^2-dn^2}\f1{m+n\sqrt d}\eq0\pmod{p^2}.
\end{equation}
\end{theorem}
\begin{remark} Our computation indicates that in general we cannot replace the modulus $p^2$
in \eqref{p2} by higher powers of $p$.
\end{remark}

We are going to prove Theorems \ref{Th1.1} and \ref{Th1.2} in Sections 2 and 3, respectively.

\section{Proof of Theorem \ref{Th1.1}}
\setcounter{lemma}{0}
\setcounter{theorem}{0}
\setcounter{equation}{0}
\setcounter{conjecture}{0}
\setcounter{remark}{0}
\setcounter{corollary}{0}

We need the following two well known results (see, e.g., \cite{IR}).

\begin{lemma}\label{Lem2.1} For any odd prime $p$, we have
\begin{equation} x^{p-1}-1\eq\prod_{n=1}^{p-1}(x-n)\pmod p.
\end{equation}
\end{lemma}

\begin{lemma}\label{Lem2.2} Let $p$ be an odd prime. For any algebraic integer $\al$, we have
the congruence
$$(x-\al)^p\eq x^p-\al^p\pmod p$$
in the ring of all algebraic integers.
\end{lemma}

\medskip
\noindent{\bf Proof of Theorem \ref{Th1.1}}. Clearly,
$1/{\sqrt d}=\sqrt d/d$. In view of Lemma \ref{Lem2.1} and the well known congruence
$$d^{(p-1)/2}\eq\l(\f dp\r)\pmod p,$$
we have
\begin{align*}\prod_{m,n=1}^{p-1}(x-m-n\sqrt d)&=\sqrt d^{(p-1)^2}\prod_{m,n=1}^{p-1}\l(\f{x-m}{\sqrt d}-n\r)=\l(d ^{(p-1)/2}\r)^{p-1}\prod_{m=1}^{p-1}\l(\l(\f{x-m}{\sqrt d}\r)^{p-1}-1\r)
\\&\eq \prod_{m=1}^{p-1}\l(\l(\f dp\r)(x-m)^{p-1}-1\r)
=\f{\prod_{m=1}^{p-1}((x-m)^p-(\f dp)(x-m))}{\prod_{m=1}^{p-1}(x-m)}
\\&\eq\f{\prod_{m=1}^{p-1}(x^p-m^p-(\f dp)(x-m))}{x^{p-1}-1}
\\&\eq\f1{x^{p-1}-1}\prod_{m=1}^{p-1}\l(x^p-\l(\f dp\r)x-\l(1-\l(\f dp\r)\r)m\r)\pmod p
\end{align*}
and hence
\begin{equation}\label{mn}\begin{aligned}\prod_{m,n=1}^{p-1}(x-m-n\sqrt d)&\eq\begin{cases}(x^p-x)^{p-1}/(x^{p-1}-1)=x^{p-1}(x^{p-1}-1)^{p-2}\pmod p&\t{if}\ (\f dp)=1,
\\((x^p+x)^{p-1}-1)/(x^{p-1}-1)\pmod p&\t{if}\ (\f dp)=-1.\end{cases}
\end{aligned}\end{equation}

 {\it Case} 1. $(\f dp)=1$.

 In this case, there is a unique number $r\in\{1,\ldots,(p-1)/2\}$ such that $r^2\eq d\pmod p$.
 For $m,n\in\{1,\ldots,p-1\}$, we have
 $$m^2-dn^2\eq0\pmod p \iff  \f mn\eq\pm r\pmod p\iff m\in\{\{nr\}_p,p-\{nr\}_p\},$$
 where we use $\{a\}_p$ to denote the least nonnegative residue of an integer $a$ modulo $p$.
 Thus
 \begin{align*}\prod_{m,n=1\atop p\mid m^2-dn^2}^{p-1}(x-m-n\sqrt d)&=\prod_{n=1}^{p-1}(x-\{nr\}_p-n\sqrt d)
 (x-(p-\{nr\}_p)-n\sqrt d)
 \\&=\prod_{n=1}^{p-1}\l((x-n\sqrt d)^2-(nr)^2\r)
 \\&\eq\prod_{n=1}^{p-1}x(x-2n\sqrt d)=(\sqrt d x)^{p-1}\prod_{n=1}^{p-1}\l(\f x{\sqrt d}-2n\r)
 \\&\eq (\sqrt d x)^{p-1}\l(\l(\f x{\sqrt d}\r)^{p-1}-1\r)=x^{p-1}\l(x^{p-1}-d^{(p-1)/2}\r)
 \\&\eq x^{p-1}(x^{p-1}-1)\pmod p.
 \end{align*}
 On the other hand, by \eqref{mn} we have
 $$\prod_{m,n=1}^{p-1}(x-m-n\sqrt d)\eq x^{p-1}(x^{p-1}-1)^{p-2}\pmod p.$$
 Therefore,
 $$P(x)=\f{\prod_{m,n=1}^{p-1}(x-m-n\sqrt d)}{\prod_{1\ls m,n\ls p-1\atop p\mid m^2-dn^2}
 (x-m-n\sqrt d)}
 \eq \f{x^{p-1}(x^{p-1}-1)^{p-2}}{x^{p-1}(x^{p-1}-1)}=(x^{p-1}-1)^{p-3}\pmod p.$$
Note that
\begin{align*}(x^{p-1}-1)^{p-3}&=\sum_{k=1}^{p-2}\bi{p-3}{k-1}x^{(k-1)(p-1)}(-1)^{p-3-(k-1)}
\\&=\sum_{k=1}^{p-2}\f{k(k+1)}{(p-1)(p-2)}\bi{p-1}{k+1}x^{(k-1)(p-1)}(-1)^{k-1}
\\&\eq\sum_{k=1}^{p-2}\f{k(k+1)}2x^{(k-1)(p-1)}\pmod p
\end{align*}
since $\bi{p-1}j\eq(-1)^j\pmod p$ for all $j=0,\ldots,p-1$.

{\it Case} 2. $(\f dp)=-1$.

In this case, $m^2-dn^2\not\eq0\pmod p$ for all $m,n\in\{1,\ldots,p-1\}$. Thus, by \eqref{mn} we have
$$P(x)\eq \f{(x^p+x)^{p-1}-1}{x^{p-1}-1}=\f{x^{p-1}(x^{p-1}+1)^p-(x^{p-1}+1)}{x^{2(p-1)}-1}\pmod p$$
and hence
$$P(x)\eq\f{x^{p-1}(x^{p(p-1)}+1)-x^{p-1}-1}{x^{2(p-1)}-1}=\f{x^{p^2-1}-1}{x^{2(p-1)}-1}
=\sum_{k=0}^{(p-1)/2}x^{2k(p-1)}\pmod p.$$

In view of the above, we have completed the proof of Theorem \ref{Th1.1}. \qed

\section{Proof of Theorem \ref{Th1.2}}
\setcounter{lemma}{0}
\setcounter{theorem}{0}
\setcounter{equation}{0}
\setcounter{conjecture}{0}
\setcounter{remark}{0}
\setcounter{corollary}{0}

We need the following well known result as a lemma (cf. \cite[IR]).

\begin{lemma} \label{Lem3.1} Let $p$ be a prime, and $k\in\N=\{0,1,2,\ldots\}$. Then
$$\sum_{n=1}^{p-1}n^k\eq\begin{cases}-1\pmod p&\t{if}\ p-1\mid k,
\\0\pmod p&\t{if}\ p-1\nmid k.\end{cases}$$
\end{lemma}

\medskip
\noindent{\bf Proof of Theorem 1.2}. As we can easily check the desired result for $p=3$, below we assume that $p\gs5$.

 For $m,n\in\{1,\ldots,p-1\}$, define
\begin{align*}
S(m,n)=&\ \frac{1}{m+n\sqrt{d}}+\frac{1}{p-m+n\sqrt{d}}+\frac{1}{m+(p-n)\sqrt{d}}
+\frac{1}{p-m+(p-n)\sqrt{d}}
\\&\ +\frac{1}{n+m\sqrt{d}}+\frac{1}{p-n+m\sqrt{d}}+\frac{1}{n+(p-m)\sqrt{d}}+\frac{1}{p-n+(p-m)\sqrt{d}}
\end{align*}
It is easy to verify that $S(m,n)=a(m,n)/b(m,n)$ with
\begin{equation}
	b(m,n)\eq (d m^2 - n^2)^2 (-m^2 + dn^2)^2 \pmod p
\end{equation}
and
\begin{equation}a(m,n)\eq p(1+\sqrt{d})f(d,m,n) \pmod {p^2},
\end{equation}
where
\begin{align*} f(d,m,n)=&\ -2dm^6- 2d^2m^6- 2m^4n^2+ 4dm^4n^2+ 4d^2m^4n^2- 2d^3m^4n^2- 2m^2n^4
\\&\ + 4dm^2n^4+ 4d^2m^2n^4-2d^3m^2n^4-2dn^6- 2d^2n^6.
\end{align*}
Therefore,
\begin{equation}
	S(m,n)= \frac{a(m,n)}{b(m,n)}\eq \frac{a(m,n)}{(dm^2-n^2)^2(-m^2+dn^2)^2} \pmod {p^2}.
\end{equation}

For integers $j,k\gs0$ with $j+k=6$, we have
\begin{align*}
	&\sum_{1\ls m,n\ls p-1\atop p\nmid (m^2-dn^2)(n^2-dm^2)}\frac{m^jn^k}{(dm^2-n^2)^2(-m^2+dn^2)^2}\\
	\eq &\sum_{1\ls m,n\ls p-1\atop p\nmid (m^2-dn^2)(n^2-dm^2)} m^jn^k(dm^2-n^2)^{p-3}(-m^2+dn^2)^{p-3}\\
	\eq& \sum_{1\ls m,n\ls p-1}m^jn^k(dm^2-n^2)^{p-3}(-m^2+dn^2)^{p-3}\pmod {p}.
\end{align*}

When we expand $((dm^2-n^2)(dn^2-m^2))^{p-3}$, we get terms of the form
$ m^sn^t$ with $s+t=4p-6$.
It is impossible that $p-1$ divides both $s$ and $t$ unless $p=3$.
If $p-1\nmid s$, then $\sum_{m=1}^{p-1}m^s\eq0\pmod p$
by Lemma \ref{Lem3.1}. If $p-1\nmid t$, then $ \sum_{n=1}^{p-1}n^t\eq 0 \pmod {p}$ by Lemma \ref{Lem3.1}.
Thus, for $s,t\in\N$ with $s+t=4p-6$, we have
$$ \sum_{1\ls m,n\ls p-1}m^sn^t\eq 0 \pmod {p}.$$
Hence
$$\sum_{1\ls m,n\ls p-1\atop p\nmid (m^2-dn^2)(n^2-dm^2)}a(m,n)(dm^2-n^2)^{p-3}(-m^2+dn^2)^{p-3}\eq 0\pmod {p^2} .$$
Therefore,
$$ \sum_{1\ls m,n\ls p-1\atop p\nmid (m^2-dn^2)(n^2-dm^2)}\frac{a(m,n)}{(dm^2-n^2)^2(-m^2+dn^2)^2}\eq 0\pmod {p^2}$$
and
\begin{align*}
	\sum_{1\ls m,n\ls p-1\atop p\nmid (m^2-dn^2)(n^2-dm^2)}\frac{1}{m+n\sqrt{d}}
	\eq&\  \frac{1}{8}\sum_{1\ls m,n\ls p-1\atop p\nmid (m^2-dn^2)(n^2-dm^2)}S(m,n)\\
	\eq&\ \frac{1}{8}\sum_{1\ls m,n\ls p-1\atop p\nmid (m^2-dn^2)(n^2-dm^2)}\frac{a(m,n)}{(d m^2 - n^2)^2 (-m^2 + dn^2)^2}\\\eq&\ 0\pmod {p^2}.
\end{align*}
This concludes our proof of Theorem \ref{Th1.2}. \qed

\medskip
\noindent{\bf Declaration of Interests}. There are no competing interests to declare.


\begin{thebibliography}{99}

\bibitem{IR} K. Ireland and M. Rosen, A Classical Introduction to Modern Number Theory,
2nd Edition, Grad. Texts. Math., vol. 84, Springer, New York, 1990.

\bibitem{K} N. Kalinin, {\it Wolstenholme's theorem over Gaussian integers},  arXiv:2504.07978, 2025.

\bibitem{W} J. Wolstenholme, {\it On certain properties of prime numbers}, Quart. J. Math. {\bf 5}
 (1862), 35--39.

\end{thebibliography}
\end{document}